\documentclass[11pt,twoside]{amsart}
\usepackage[all]{xy}
        \usepackage {amssymb,latexsym,amsthm,amsmath,mathtools}
        \usepackage{enumitem,color}
        
\long\def\symbolfootnote[#1]#2{\begingroup%
\def\thefootnote{\fnsymbol{footnote}}\footnote[#1]{#2}\endgroup}
\newcommand{\Z}{\ensuremath{\mathcal{Z}}}

\newcommand{\tra}{\ensuremath{{}^t}}

\makeatletter
\def\imod#1{\allowbreak\mkern10mu({\operator@font mod}\,\,#1)}
\makeatother

\newtheorem{theorem}{Theorem}[section]

\newtheorem*{theorem*}{Theorem}
\theoremstyle{definition}
\newtheorem{definition}[theorem]{Definition}
\newtheorem{problem}[theorem]{Problem}

\newtheorem{example}[theorem]{Example}

\numberwithin{equation}{section}

\newcommand{\ignore}[1]{}

\newcommand{\mynote}[1]{}
\def \GL {\mathrm{GL}}
\def \SL {\mathrm{SL}}
\def \Sp {\mathrm{Sp}}
\def \O {\mathrm{O}}
\def \U {\mathrm{U}}
\def \PSp {\mathrm{PSp}}

\begin{document}
\setcounter{section}{0}
\title{$z$-classes in groups: a survey}
\author{Sushil Bhunia}
\email{sushilbhunia@gmail.com}
\address{IISER Mohali, Knowledge City,  Sector 81, S.A.S. Nagar, Punjab 140306, India}
\author{Anupam Singh}
\email{anupamk18@gmail.com}
\address{IISER Pune, Dr. Homi Bhabha Road, Pashan, Pune 411 008 India}
\thanks{The second named author would like to acknowledge support of SERB core research grant CRG/2019/000271 for this work. This research was supported in part by the International Centre for Theoretical Sciences (ICTS) during a visit for participating in the program -  Group Algebras, Representations and Computation ICTS/Prog-garc2019/10.}
\today

\dedicatory{Dedicated to Professor I. B. S. Passi on the occasion of his $80^{th}$ birthday}

\begin{abstract}
This survey article explores the notion of $z$-classes in groups. The concept introduced here  is related to the notion of orbit types in transformation groups, and types or genus in the representation theory of finite groups of Lie type. Two elements in a group are said to be $z$-equivalent (or $z$-conjugate) if their centralizers are conjugate. This is a weaker notion than the conjugacy of elements. In this survey article, we present several known results on this topic and suggest some further questions.
\end{abstract}

\maketitle
 
\section{Introduction}

Let $G$ be a group. Two elements $g_1$ and $g_2$ are called {\bf conjugate} in $G$ if there exists a $t$ in $G$ such that $tg_1t^{-1}=g_2$. Conjugacy is an equivalence relation which gives rise to the conjugacy classes. For a finite group $G$, the number of conjugacy classes is same as the number of non-equivalent irreducible complex-representations. Thus, computing conjugacy classes is one of the central problems in group theory and representation theory. Let $g$ be an element of $G$. The centralizer of $g$ in $G$ is denoted as 
$$\mathcal Z_G(g)=\{x\in G\mid xg=gx\}.$$ 
For a group $G$, we say, two elements $g_1$ and $g_2$ are {\bf $z$-equivalent or $z$-conjugate} if their centralizers are conjugate subgroups within $G$, i.e., if there exists a $t$ in $G$ such that $t\mathcal Z_G(g_1)t^{-1} = \mathcal Z_G(g_2)$. Clearly, $z$-equivalency is an equivalence relation on $G$, and represents the conjugacy classes of centralizer subgroups of $G$. These equivalence classes are called $z$-classes or centralizer classes. It is easy to verify that if two elements are conjugate then they are $z$-equivalent. However, the converse need not be true. 
More precisely, $z$-classes are union of conjugacy classes. The notion of $z$-classes was introduced by Ravi Kulkarni (see~\cite{Ku, Ku1}). This notion is different from asking if the centralizers are isomorphic (abstractly) groups. We begin with some examples to understand the complications involved. 
\begin{example}
For an Abelian group, there are as many conjugacy classes as the number of elements. However, there is only one $z$-class. 
\end{example}
\noindent Thus, we focus on studying non-Abelian groups only.
\begin{example}\label{example1.2}
For the symmetric group $S_3$ and $S_4$, the $z$-classes are same in number as the number of conjugacy classes. We know that for the symmetric group $S_n$, the conjugacy classes correspond to the partitions of $n$. However, in $S_5$, the elements $(123)(45)$ and $(123)$ are $z$-conjugate, as they have the same centralizer, but not conjugate. The precise result for the symmetric groups $S_n$ and alternating group $A_n$ is in~\cite{BKS} and we discuss this briefly in Section~\ref{z-classes-sym}.
\end{example}
\begin{example}\label{GLn}
In the group $\GL_n(k)$, the diagonal matrices with distinct diagonal entries (these are example of regular semisimple elements) have centralizer the diagonal maximal torus, thus they are all $z$-equivalent even though they are not conjugate. More generally, $z$-classes of regular semisimple elements is same as the conjugacy classes of maximal tori. The $z$-classes for $\GL_n(\mathbb F_q)$ was studied by Green in~\cite{Gr} in the context of theory of types while computing characters for this group. More general case of matrices over division rings and affine linear transformations have been studied in~\cite{Go} and~\cite{Ku}, respectively. The notion of types is extended to ``generalised types'' for the set of matrices $\mathrm{M}_n(k)$ in Britnell and Wildon~\cite{BW1,BW2} which coincides with the notion of $z$-classes (extended to an associative algebra) in this case.
\end{example}
\begin{example}
In the quaternion group $Q_8=\{\pm 1, \pm i, \pm j, \pm k\}$, the centralizers of $i$ and $j$ are abstractly isomorphic as groups but not conjugate. 
\end{example} 

Thus, the main problem here is to list the representatives of $z$-classes for a given group, and possibly compare it with the conjugacy classes. Our focus will be mainly on linear groups (the groups which are subgroups of matrix groups). In Section~\ref{section-GM}, we begin with the notion of {\bf orbit types} in geometry and mention the known results about $z$-classes for some Lie groups.  In Section~\ref{section-CG}, we discuss the classical groups and the finiteness results over the fields of type (F). These results were proved case-by-case with explicit computation of centralizers. In Section~\ref{section-AG}, we list the finiteness result known for reductive algebraic groups over fields of type (F) using Galois cohomology. Finally, in Section~\ref{section-FG}, we put together the results known for some other finite groups, especially $p$-groups and Weyl groups. In this survey article, we have summarised the results known on this problem as far as we could gather information. All throughout we suggest problems which, as far as we know, is still to be solved. We hope this article will be useful to the graduate students and researchers in this area.

\subsection*{Acknowledgement:}
The second named author gratefully acknowledges the opportunity to attend a series of lectures in his graduate days given at HRI by Ravi Kulkarni, introducing the notion of $z$-classes.  

\section{Geometrical motivation and $z$-classes for Lie groups}\label{section-GM}

Let $G$ be a compact group, and $X$ be a $G$-space. Two orbits $G.x$ and $G.y$ have the {\bf same orbit type} if the stabilisers $G_x$ and $G_y$ are conjugate in $G$ (see~\cite[Definition 10.5]{HM}). We can use this to define an equivalence relation on $X$ by saying that; $x,y\in X$ have same orbit types if $G_x$ and $G_y$ are conjugate in $G$. The main question is to determine the orbit types (representatives of the orbit types). This is used to define stable isotropy, that is, $X$ has a single orbit type. This concept is a generalisation of free action. Montgomery (see~\cite[Problem 45]{Ei}) conjectured that if $G$ is a compact Lie group acting on a compact manifold, then it has only finitely many orbit types. Floyd (see~\cite{Flo}) proved this conjecture positively for a torus group acting on a compact orientable manifold, and Mostow~\cite{mos} proved the conjecture fully. A further generalisation of this result is due to Mann~\cite{Ma}. We give an example here.
\begin{example}
Consider the compact real group $\mathrm{SO}(n)$ acting on $(n-1)$-sphere $\mathbb S^{n-1}$, for $n\geq 3$, as follows: $\mathrm{SO}(n)\times \mathbb S^{n-1} \rightarrow \mathbb S^{n-1}$ given by $(g, x) \mapsto g(x)$. Then, the isotropy groups are conjugate to a subgroup $\mathrm{SO}(n-1)$.
\end{example}

The study of $z$-classes  can be thought of as a special case where $G$ acts on itself by conjugation. Thus, from the results mentioned above, it follows that the number of $z$-classes is finite for a compact Lie group. Now, a natural question arises in this context is that, what is the (precise) number of $z$-classes? In an attempt to answer this question, Singh~\cite{Si} computed the number of $z$-classes for the real compact group of type $G_2$. Bose~\cite{Bo} calculated explicitly the number of $z$-classes for compact simple Lie groups of type $A_n$, $B_n$, $C_n$, $D_n$, $F_4$ and $G_2$ case-by-case. For details, we refer the reader to~\cite[Table in Section 8]{Bo}. It would be nice to complete this computation for the remaining simple Lie groups. Thus we suggest: 
\begin{problem}
Compute the number of $z$-classes for the compact Lie groups of type $E_6, E_7, E_8$.
\end{problem}

The rank one symmetric spaces of compact types are: sphere, complex and quaternionic projective spaces, and octonionic projective plane. The corresponding isometry groups are the following real compact Lie groups: $\O(n+1)$, $\mathrm{PSU}(n+1)\rtimes \mathbb{Z}/2\mathbb Z$, $\PSp(n+1)$ and compact group of type $F_4$. Since, these groups are compact Lie groups of type: $A_n$, $B_n$, $C_n$, $D_n$ and $F_4$, from the work of Bose mentioned above, we already have the precise number of $z$-classes in these cases. Now, the rank one symmetric spaces of non-compact types are: real, complex and quaternionic hyperbolic spaces, and the octonionic hyperbolic plane. In this case, the corresponding isometry groups are the following signature one groups: $\mathrm{PO}(n,1)$, $\mathrm{PU}(n,1)$, $\mathrm{PSp}(n,1)$ and the group of type $F_{4(-20)}$. Gongopadhyay and Kulkarni~\cite{GK1} studied $z$-classes of $\mathrm{PO}(n,1)$, the isometry group of the real hyperbolic space. Gongopadhyay~\cite[Theorem 1.1]{Go1} gave an explicit counting of $z$-classes in $\Sp(n,1)$, the isometry group of the quaternionic hyperbolic space. For an explicit counting of $z$-classes in $\mathrm{U}(n,1)$, the isometry group of complex hyperbolic space, we refer an interested reader to~\cite[Proposition 4.4]{BS}. Thus, it remains to do the following.
\begin{problem}
Compute the number of $z$-classes for the group of type $F_{4(-20)}$, the isometry group of  octonionic hyperbolic plane.
\end{problem}

\section{$z$-classes for classical groups}\label{section-CG}

In this section, we explore the $z$-classes for classical groups. For the definition of these groups, we refer a reader to the book by Grove~\cite{Gv} on the subject. We begin with an example to highlight a particular aspect of our problem. 
\begin{example}
From the theory of Jordan canonical forms, it is clear that the group $\GL_2(\mathbb C)$ has only finitely $z$-classes. In fact, it has exactly three $z$-classes given by the representatives $I, \begin{pmatrix} 1&0\\0&2\end{pmatrix}, \begin{pmatrix} 1&1\\0&1\end{pmatrix}$. However, the group $\GL_2(\mathbb Q)$ has infinitely many $z$-classes of regular semisimple elements caused by the presence of infinitely many non-conjugate maximal tori in this group.
\end{example}
\noindent This example leads us to believe that the arithmetic nature of the base field plays an important  role in finiteness of $z$-classes. The theorem of Steinberg (which we mention in the following section when we work with algebraic groups, see Theorem~\ref{Steinberg}) tells us so over an algebraically closed field. We recall the definition of {\bf fields of type (F)} due to Borel and Serre~\cite{BSe}. 
\begin{definition}[Fields of type (F)]\label{fields-type-F}
A perfect field $k$ is said to be of type (F) (or said to have the property (F)) if $k$ has only finitely many field extensions of any finite degree, up to isomorphism. 
\end{definition}
\noindent Examples of such fields are, algebraically closed fields (for example, $\mathbb C$), real numbers $\mathbb R$, local fields (for example, $\mathbb Q_p$), and finite fields $\mathbb F_q$.  The field $\mathbb Q$ or any number field does not have the property (F). As we mentioned in the introduction, Kulkarni~\cite[Theorem 7.4]{Ku}, proved that the number of $z$-classes in $\GL_n(k)$ is finite when $k$ is of type (F). 

Let $V$ be an $n$-dimensional vector space over a field $k$ of characteristic $\neq 2$ which is of type (F), equipped with a non-degenerate symmetric or skew-symmetric bilinear form $B$. Gongopadhyay and Kulkarni~\cite[Theorem 1.1]{GK} proved that the number of $z$-classes in orthogonal groups $\O(V, B)$ and symplectic groups $\Sp(V, B)$ are finite. 

Suppose, $k$ is a perfect field of characteristic $\neq 2$ with a non-trivial Galois automorphism of order $2$. Let $V$ be a finite dimensional vector space over $k$ with a non-degenerate Hermitian form $B$. Suppose, the fixed field $k_0$ is of type (F). Then, in~\cite[Theorem 1.2]{BS}, we proved that the number of $z$-classes in the unitary group $\U(V, B)$ is finite. The condition that field is of type (F) is necessary for the above results as shown by the following example (see~\cite[Example 5.1]{BS}).
\begin{example}
We embed $\GL_2(\mathbb Q)$ in $\Sp_4(\mathbb Q)$ with respect to the skew-symmetric form $\begin{pmatrix} & I_2 \\ -I_2 & \end{pmatrix}$ given by $A\mapsto \begin{pmatrix} A & \\ & \tra {A}^{-1}\end{pmatrix}$. This embedding gives maximal tori in $\Sp_4(\mathbb Q)$ coming from that of $\GL_2(\mathbb Q)$. This leads to infinitely many $z$-classes (of semisimple elements) in $\Sp_4(\mathbb Q)$. A similar example could be constructed for orthogonal and unitary groups. For details, we refer an interested reader to~\cite[Section 5]{GS}.
\end{example}
\noindent The conjugacy classes for these groups have been known, and most of the work mentioned in this section relies on explicitly describing the centralizers of elements, which is useful in its own right. The problem of studying $z$-classes in classical groups have two main components. First, whether the number of $z$-classes is finite or not. The second component is to compute the exact number of $z$-classes, whenever the number is finite or the group in question is finite. So far, we have shown that the classical groups over fields of type (F) have finitely many $z$-classes. Now we move on to the counting aspect of this problem. 

\subsection{Counting and types}

If we look at the character table of $\SL_2(q)$ (one can find this in~\cite{B} or~\cite{Pr}) we notice that several conjugacy classes and irreducible characters are grouped together in parametrized form. One observes a similar pattern in $\Sp_4(q)$ from the work of Srinivasan~\cite{Sr}. In~\cite{Gr}, Green studied the complex representations of $\GL_n(q)$ where he introduced the function $t(n)$ for the `types of characters/classes' which is the number of $z$-classes in $\GL_n(q)$. In Deligne-Lusztig theory, where we study the representation theory of finite groups of Lie type, $z$-classes of semisimple elements play an important role. Carter~\cite{Ca} and Humphreys~\cite{Hu2} define genus of an element in algebraic group $G$ over $k$. This definition was generalised by Bose in~\cite{Bo}. Two elements have the same genus if they are $z$-equivalent in $G(k)$ and the genus number (respectively semisimple genus number) is the number of $z$-classes (respectively the number of $z$-classes of semisimple elements). While calculating genus number for simply connected simple algebraic groups over algebraically closed field (see \cite[Table in Section 8]{Bo}), Bose computed $z$-classes for some of these groups too. Thus, understanding $z$-classes for finite groups of Lie type, especially semisimple genus, and their counting is of importance in representation theory as well (see~\cite{Fl, FG, FS, Ca, DM}). We finish this section by recalling a result which fits in this discussion.
\begin{theorem}\cite[Theorem 1.3]{BS}
For $q>n$, the number of $z$-classes in $\U_n(q)$ is same as the number of $z$-classes in $\GL_n(q)$. Further, the generating function for the number of $z$-classes is 
$\displaystyle\prod_{i=1}^{\infty} z(x^i)$, where $z(x)=\displaystyle\prod_{j=1}^{\infty}\frac{1}{(1-x^j)^{p(j)}}$ 
and $p(j)$ is the number of partitions of $j$.
\end{theorem}
With this in mind, we suggest some questions in this direction for further exploration. 
\begin{problem}
Determine the number of $z$-classes in $\Sp_n(q)$ and $\O_n(q)$ and write the generating functions.
\end{problem}
\begin{problem}
Compute the number of $z$-classes in the finite exceptional groups of Lie type.
\end{problem}

\subsection{Classical groups over rings}

We can consider the classical groups over certain rings as well, for example, PIDs, EDs, or local rings. For definitions of the classical groups in this generality, we refer a reader to the book~\cite{HO}. The centralizers of elements for some of these groups have been studied in~\cite{Sg, PSS} over local rings of length two. We suggest the following question.
\begin{problem}
Parametrize the $z$-classes in classical groups over local rings of length two.
\end{problem}

\section{$z$-classes for algebraic groups}\label{section-AG}

We are interested in linear algebraic groups defined over a base field $k$. Let $K$ be an algebraically closed field. An {\bf algebraic group} $G$ is an affine variety defined over $K$ with a group structure where the two group operations, namely multiplication and inversion, are morphisms of varieties. This definition is equivalent to having a Hopf algebra structure on the coordinate ring $K[G]$. Let $k$ be a field and its algebraic closure is $\bar k$. An algebraic group $G$ (a priory over $\bar k$) is said to be defined over base field $k$ if the coordinate ring $\bar k[G]$ is defined over $k$. The $k$-points of $G$ are denoted as $G(k)$. For a detailed account of the theory of linear algebraic groups, one can see the classic textbooks~\cite{Hu, Sp} among many others.

Now, we explore the question of conjugacy classes and $z$-classes in $G(k)$ and some of the subtleties involved. We remark that when we look for conjugacy in this context, the conjugating elements are in $G(k)$. The structure of centralizers plays an important role in the classification of linear algebraic groups. For example, centralizer of a semisimple element in a connected semisimple group is again reductive (see Chapter 2 in~\cite{Hu2}). This helps us in proving results using the method of induction. Once again, when we consider the $z$-classes we take the centralizer subgroup and its $k$-point under the conjugation action. The issue of smoothness etc could arise over arbitrary field. We begin with an interesting result by Steinberg (Section 3.6 Corollary 1 to Theorem 2 in~\cite{St}):
\begin{theorem}\label{Steinberg}
Let $G$ be a connected reductive algebraic group defined over an algebraically closed field $K$ in good characteristic. Then the number of conjugacy classes of centralizers in $G$ (i.e., $z$-classes) is finite. 
\end{theorem} 
\noindent Notice that the number of conjugacy classes could be infinite in these cases (see Example~\ref{GLn}).  Bose calculated genus number for simply connected simple algebraic groups over algebraically closed field (see \cite[Table in Section 8]{Bo}), thus completing the counting problem here.
We would like to raise the following question:
\begin{problem}
Let $G$ be a linear algebraic group defined over a field $k$. Classify the $z$-classes in $G(k)$.
\end{problem} 
\noindent In~\cite{GS}, this problem is studied for reductive groups and the following is proved with some mild assumptions on $k$ and $G$ to ensure smoothness of centralizers and normalizers (see~\cite[Theorem 4.3]{GS}).
\begin{theorem}
 Let $G$ be a reductive algebraic group defined over a field $k$ of type (F). Then, the group $G(k)$ has finitely many $z$-classes. 
\end{theorem}
\noindent Now the problem remains to show if we can extend the above result to arbitrary algebraic groups, a particular case being solvable groups and nilpotent groups.

\subsection{$z$-classes for solvable groups}

All the results stated so far in this section are for reductive algebraic groups. Now, we move on to discuss other extremes, namely solvable linear algebraic groups, which is not much explored yet. 
Let $\mathrm{B}_n(k)$ denote the group of upper triangular matrices in $\mathrm{GL}_n(k)$.
Lie-Kolchin theorem says that a connected solvable linear algebraic group $G$ defined over an algebraically closed field $k$ is a subgroup of $\mathrm{B}_n(k)$ for some $n$. Before we state the main result of this section, we begin with a simple example. For the remainder of this section, we assume that $k$ is an arbitrary field of characteristic $\neq 2$.
\begin{example}
We know that the number of unipotent conjugacy classes in $\GL_n(k)$ is in one to one correspondence with the number of partitions of $n$. However, in general, this is not true for its subgroups. For instance, $\mathrm{B}_3(k)$ the subgroup of upper triangular matrices in $\GL_3(k)$ has $5$ unipotent conjugacy classes. Note that, the number of semisimple $z$-classes in $\mathrm{B}_3(k)$ is $5$. Thus, the number of $z$-classes in $\mathrm{B}_3(k)$ is finite (for details, see \cite[Proposition 2.1 and Appendix]{Bh}). 
\end{example}
In general, this gives a misleading picture as we will see soon. Even the finiteness result, as for the reductive algebraic groups, depending on the arithmetic nature of the base field doesn't hold true. More precisely, it turns out that the number of $z$-classes in $\mathrm{B}_n(k)$ is infinite whenever the field $k$ is infinite and $n\geq 6$ (see~\cite[Theorem 1.2]{Bh}). It is finite otherwise. Thus a more general question remains:
\begin{problem}
For a nilpotent and solvable linear algebraic group defined over $k$, determine when the number of $z$-classes is finite (possibly depending on derived length).
\end{problem}

\section{$z$-classes for finite groups}\label{section-FG}
In this section, unless specified otherwise, we assume that the groups are finite.  We begin with some examples.
\begin{example}
Let $D_{2n}=\langle r, s\mid r^n=s^2=1, rs=sr^{-1}\rangle$ be the dihedral group of order $2n$. The element $r$ and all its powers are called rotations and $r^is$ are called reflections.
When $n$ is odd, the number of $z$-classes in $D_{2n}$ is $3$, corresponding to the identity element, $r$ and $s$. When $n\equiv 2 \imod 4$, the $\mathcal{Z}_{D_{2n}}(s)=\mathcal{Z}_{D_{2n}}(sr^{n/2})$. Thus, two non-conjugate reflections are $z$-conjugate. Hence the number of $z$-classes is $3$ again, corresponding to the identity, $r$ and $s$, in this case as well. But when $n\equiv 0\mod 4$, the two distinct conjugacy classes of reflections are not $z$-conjugate. Thus, the number of $z$-classes in $D_{2n}$ is $4$ corresponding to central elements, $r$ and two reflections, in this case. Therefore, we have the following: 
\begin{equation*}
\text{Number of $z$-classes in $D_{2n}$}=\left\{\begin{array}{ll}
3 & \text{ if } n\equiv 1, 2, 3\mod 4,\\
4 & \text{ if } n \equiv 0\mod 4.
\end{array}\right.
\end{equation*}
\end{example}

Clearly, if two groups are isomorphic then they have the same number of $z$-classes. In fact, it is true under certain weaker assumptions. Let us recall the definition.
\begin{definition}
Two groups $G$ and $H$ are said to be \emph{isoclinic} if there is an isomorphism $\varphi\colon G/\Z(G)\rightarrow H/\Z(H)$ and an isomorphism $\psi \colon [G,G]\rightarrow [H,H]$ such that the following diagram commutes:
\[\xymatrix{
G/\Z(G)\times G/\Z(G) \ar[r]^{\varphi\times \varphi} \ar[d]_{\pi_G}  & H/\Z(H)\times H/\Z(H)\ar[d]^{\pi_H} \\ [G,G] \ar[r]^{\psi} &[H,H]},\]
where the vertical maps $\pi_G$ and $\pi_H$ are the natural commutator maps.
\end{definition}
Kulkarni et al., in \cite[Theorem 2.2]{KKJ}, proved that if $G$ and $H$ are isoclinic, then the number of $z$-classes in $G$ is equal to the number of $z$-classes in $H$. We know that a normal subgroup is a union of conjugacy classes. Likewise, in any group $G$, a maximal Abelian normal subgroup is a union of $z$-classes (see \cite[Theorem 3.5]{KKJ}).

On this note, we also point out that the classification of finite groups with some given condition on the set of centralizers is studied. Let $\mathcal C(G)$ be the set of centralizers of $G$. In~\cite{Sz, DHJ,Za} and~\cite{KZ} the classification of groups is studied with given condition on $|\mathcal C(G)|$.

A very important family of finite groups is $p$-groups, which we will deal with in the subsection to follow.

\subsection{$z$-classes for finite $p$-groups}

Here we collect some results regarding study of finite $p$-groups using the $z$-classes. Kulkarni et al., in \cite[Theorem 3.7, 3.10]{KKJ}, proved that for a non-Abelian finite $p$-group $G$ the number of $z$-classes is $p+2$ if and only if either $G/\Z(G)\cong C_p\times C_p$ or $G$ has a unique Abelian subgroup of index $p$ and the order of $\Z(G/\Z(G))$ is $p$, where $C_p$ is a cyclic group of order $p$. If $G$ is a group of order $p$ or $p^2$ then $G$ is Abelian, thus the number of $z$ classes is $1$. Now, if $G$ is a non-Abelian group of order $p^3$ or $p^4$ then the number of $z$-classes is $p+2$. Jadhav and Kitture in~\cite{JK} studied the $z$-classes in groups of order $p^5$.
Thus, the following question remains:
\begin{problem}
 Determine the $z$-classes in groups of order $p^n$ for $n\geq 6$.
\end{problem}

 In general, Kulkarni et al. (see \cite[Theorem 3.14]{KKJ}) proved that if a finite $p$-group $G$ with $[G:\Z(G)]=p^k\; (k\geq 2)$ having $\frac{p^k-1}{p-1}+1$ many $z$-classes, then either $G/\Z(G)\cong C_p\times C_p$ or $G$ is isoclinic to special $p$-group with no Abelian subgroup of index $p$ . Observe that this is a necessary condition but not sufficient (for details, see~\cite[Remark 3.15]{KKJ}). Later, Arora and Gongopadhyay in~\cite[Theorem 1.1]{AG}, gave a sufficient condition to the above theorem when $G$ is a non-Abelian finite group of conjugate type $(n,1)$. For the sake of clarification, here we include the definition of {\bf conjugate type}. 
 \begin{definition}
A finite group $G$ is called a group of  conjugate type $(n_1, n_2, \ldots, n_k)$ if $n_1>n_2>\cdots>n_k=1$ are the indices of the centralizers of elements of $G$ in $G$.
 \end{definition}
\noindent Thus, we can ask the following more general question:
\begin{problem}
 Characterize the $z$-classes in groups of conjugate type $(n_1, n_2, \ldots, n_k)$.
\end{problem}
 
\subsection{$z$-classes for symmetric and alternating groups}\label{z-classes-sym}
It is well-known that the conjugacy classes in symmetric group $S_n$ is given by the partitions of $n$. Using combinatorial arguments related to the partitions of $n$ and centralizers structure, in \cite[Corollary 1.2]{BKS}, we computed the number of $z$-classes in $S_n$, which is equal to 
$$p(n)-p(n-2) + p(n-3) + p(n-4) - p(n-5),$$ 
where $p(m)$ is the number of partitions of $m$. From Example~\ref{example1.2}, note that the number of conjugacy classes and $z$-classes in $S_3$ and $S_4$ are same, and the number is $3$ and $5$ respectively. There is no other coincidence for higher $n$. The $z$-classes in the alternating group $A_n$ turns out to be slightly tricky, but the counting is done in a similar fashion. It follows from this work that the number of $z$-classes in $A_n$ is 
$$\mathrm{cl}(A_n)-(q(n)+\widetilde q(n-3)) + \epsilon(n),$$
where $\mathrm{cl}(A_n)=\frac{p(n)+3q(n)}{2}$ is the number of conjugacy classes in $A_n$, $q(n)$ is the number of partitions of $n$ which has all parts distinct and odd, $\widetilde q(m)$ is the number of restricted partitions of $m$, with all parts distinct, odd and which do not have $1$ (and $2$) as its part and $\epsilon(n)$ is the number of partitions of $n$ with all of its parts distinct, odd and square (for details, see \cite[Corollary 1.4]{BKS}). We suggest the following question:
\begin{problem}
Compute the number of $z$-classes for other Weyl groups. 
\end{problem}



\begin{thebibliography}{99}
\normalsize
\bibitem[AG]{AG} Arora, Shivam; Gongopadhyay, Krishnendu, \emph{``$z$-classes in finite groups of conjugate type $(n,1)$"}, Proc. Indian Acad. Sci. Math. Sci. 128 (2018), no. 3, Paper No. 31, 7 pp. 
\bibitem[Bh]{Bh} Bhunia, Sushil, \emph{``Conjugacy classes of centralizers in the group of upper triangular matrices"}, J. Algebra Appl. 19 (2020), no. 1, 2050008, 14 pp.
\bibitem[B]{B} Bonnaf\'{e}, C\'{e}dric, \emph{``Representations of $\mathrm{SL}_2(\mathbb F_q)$"}, Algebra and Applications, 13. Springer-Verlag London, Ltd., London, 2011.
\bibitem[BKS]{BKS} Bhunia, Sushil; Kaur, Dilpreet; Singh, Anupam, \emph{``$z$-classes and rational conjugacy classes in alternating groups"}, J. Ramanujan Math. Soc. 34 (2019) no. 2, 169-183.
\bibitem[Bo]{Bo} Bose, Anirban, \emph{``On the genus number of algebraic groups"}, J. Ramanujan Math. Soc. 28 (2013), no. 4, 443-482.
\bibitem[BS]{BS} Bhunia, Sushil; Singh, Anupam, \emph{``Conjugacy classes of centralizers in unitary groups"}, J. Group Theory, 22 (2019), 231-151.
\bibitem[BSe]{BSe} A. Borel and J.-P. Serre, \emph{``Th\'{e}or\`{e}mes de finitude en cohomologie galoisienne"}, (French) Comment. Math. Helv. 39 (1964), 111-164.
\bibitem[BW1]{BW1} Britnell, John R.; Wildon, Mark, \emph{``On types and classes of commuting matrices over finite fields"}, J. Lond. Math. Soc. (2) 83 (2011), no. 2, 470-492. 
\bibitem[BW2]{BW2} Britnell, John R.; Wildon, Mark, \emph{``On types of matrices and centralizers of matrices and permutations"}, J. Group Theory 17 (2014), no. 5, 875-887. 
\bibitem[Ca]{Ca} Carter, Roger W., \emph{``Finite groups of Lie type. Conjugacy classes and complex characters"}, Reprint of the 1985 original. Wiley Classics Library. A Wiley-Interscience Publication. John Wiley \& Sons, Ltd., Chichester, 1993.
\bibitem[DHJ]{DHJ} Dolfi, Silvio; Herzog, Marcel; Jabara, Enrico, \emph{``Finite groups whose noncentral commuting elements have centralizers of equal size''}, Bull. Aust. Math. Soc. 82 (2010), no. 2, 293-304. 
\bibitem[DM]{DM} Digne, Fran\c{c}ois; Michel, Jean, \emph{``Representations of finite groups of Lie type"}, London Mathematical Society Student Texts, 21. Cambridge University Press, Cambridge, 1991.
\bibitem[Ei]{Ei}S. Eilenberg, \emph{``On the problems of topology"}, Ann. of Math. (2)  50  (1949) 247-260.
\bibitem[Flo]{Flo} E. E. Floyd, \emph{``Orbits of torus groups operating on manifolds"}, Ann. of Math. (2)  65  (1957) 505-512.
\bibitem[Fl]{Fl} Fleischmann, Peter, \emph{``Finite fields, root systems, and orbit numbers of Chevalley groups"}, Finite Fields Appl. 3 (1997), no. 1, 33-47. 
\bibitem[FG]{FG} Fulman, Jason; Guralnick, Robert, \emph{``The number of regular semisimple conjugacy classes in the finite classical groups"}, Linear Algebra Appl. 439 (2013), no. 2, 488-503. 
\bibitem[FS]{FS} Fong, Paul; Srinivasan, Bhama, \emph{``The blocks of finite general linear and unitary groups"}, Invent. Math. 69 (1982), no. 1, 109-153. 
\bibitem[Go1]{Go1} Gongopadhyay, Krishnendu, \emph{``The z-classes of quaternionic hyperbolic isometries"}, J. Group Theory 16 (2013), 941-964.
\bibitem[GK]{GK} Gongopadhyay, Krishnendu and Kulkarni, Ravi S.,  \emph{``The $z$-classes of isometries"}, J. Indian Math. Soc. (N.S.) 81 (2014), no. 3-4, 245-258.
\bibitem[GK1]{GK1} Gongopadhyay, Krishnendu and Kulkarni, Ravi S., \emph{``z-classes of isometries of the hyperbolic space"}, Conform. Geom. Dyn. 13 (2009), 91-109.
\bibitem[Go]{Go} Gouraige, Rony, \emph{``z-classes in central simple algebras"}, Thesis (Ph.D.)-City University of New York. 2006.
\bibitem[Gr]{Gr} Green, J. A., \emph{``The characters of the finite general linear groups"}, Trans. Amer. Math. Soc. 80 (1955), 402-447. 
\bibitem[GS]{GS} Garge, Shripad M.; Singh, Anupam, \emph{``Finiteness of $z$-classes in reductive groups"}, J. Algebra 554 (2020), 41-53. 
\bibitem[Gv]{Gv} Larry C. Grove, \emph{``Classical groups and geometric algebra"}, Graduate Studies in Mathematics, 39. American Mathematical Society, Providence, RI, 2002.
\bibitem[HM]{HM} Hofmann, Karl H.; Morris, Sidney A., \emph{``The structure of compact groups.
A primer for the student-a handbook for the expert"}, Third edition, revised and augmented. De Gruyter Studies in Mathematics, 25. De Gruyter, Berlin, 2013. xxii+924 pp. 
\bibitem[HO]{HO} Hahn, Alexander J.; O'Meara, O. Timothy, \emph{``The classical groups and K-theory. With a foreword by J. Dieudonne"}, Grundlehren der Mathematischen Wissenschaften [Fundamental Principles of Mathematical Sciences], 291. Springer-Verlag, Berlin, 1989. xvi+576 pp.
\bibitem[Hu]{Hu} Humphreys, James E., \emph{``Linear algebraic groups"}, Graduate Texts in Mathematics, No. 21. Springer-Verlag, New York-Heidelberg, 1975. xiv+247 pp.
\bibitem[Hu2]{Hu2} Humphreys, James E., \emph{``Conjugacy classes in semisimple algebraic groups"}, Mathematical Surveys and Monographs, 43. American Mathematical Society, Providence, RI, 1995.
\bibitem[JK]{JK} Jadhav, Vikas S.; Kitture, Rahul Dattatraya, \emph{``z-classes in $p$-groups of order $\leq p^5$"}, Bull. Allahabad Math. Soc. 29 (Part 2) (2014) 173-194.
\bibitem[KKJ]{KKJ} Kulkarni, Ravindra; Kitture, Rahul Dattatraya; Jadhav, Vikas S., \emph{``$z$-classes in groups"}, J. Algebra Appl. 15 (2016), no. 7, 1650131.
\bibitem[Ku]{Ku} Kulkarni, Ravi S., \emph{``Dynamics of linear and affine maps"}, Asian J. Math. 12 (2008), no. 3, 321-344. 
\bibitem[Ku1]{Ku1} Kulkarni, Ravi S., \emph{``Dynamical types and conjugacy classes of centralizers in groups"}, J. Ramanujan Math. Soc. 22 (2007), no. 1, 35-56.
\bibitem[KZ]{KZ} K. Khoramshahi; M. Zarrin, \emph{``Groups with the same number of centralizers"}, to appear in the J. Algebra Appl.  https://doi.org/10.1142/S0219498821500122.
\bibitem[Ma]{Ma} Mann, L. N., \emph{``Finite orbit structure on locally compact manifolds"}, Michigan Math. J. 9 (1962), 87-92. 
\bibitem[Mo]{mos} G.D. Mostow, \emph{``On a conjecture of Montgomery"}, Ann. of Math., 65(2), 1957, 513-516.
\bibitem[PSS]{PSS} Prasad, Amritanshu; Singla, Pooja; Spallone, Steven, \emph{``Similarity of matrices over local rings of length two"}, Indiana Univ. Math. J. 64 (2015), no. 2, 471-514. 
\bibitem[Pr]{Pr} Prasad, Amritanshu, \emph{``Representations of $\GL_2(\mathbb F_q)$ and $\SL_2(\mathbb F_q)$, and some remarks about $\GL_n(\mathbb F_q)$"},  arXiv:0712.4051.
\bibitem[Sg]{Sg} Singla, Pooja, \emph{``On representations of general linear groups over principal ideal local rings of length two"}, J. Algebra 324 (2010), no. 9, 2543-2563. 
\bibitem[Sz]{Sz} Suzuki, Michio, \emph{``Finite groups with nilpotent centralizers''}, Trans. Amer. Math. Soc. 99 (1961), 425-470.
\bibitem[Sr]{Sr} Srinivasan, Bhama, \emph{``The characters of the finite symplectic group $\Sp(4,q)$"}, Trans. Amer. Math. Soc. 131 (1968) 488-525. 
\bibitem[Si]{Si} Singh, Anupam, \emph{``Conjugacy  Classes  of Centralizers  in $G_{2}$"}, J. Ramanujan Math. Soc. 23 (2008), no. 4, 327-336.
\bibitem[Sp]{Sp} Springer, T. A., \emph{``Linear algebraic groups"}, Second edition. Progress in Mathematics, 9. Birkhäuser Boston, Inc., Boston, MA, 1998. xiv+334 pp.
\bibitem[St]{St} Steinberg, Robert, \emph{``Conjugacy Classes in Algebraic Groups"}, notes by V. Deodhar, Lecture Notes in Mathematics 366, Springer-Verlag (1974).
\bibitem[Za]{Za} Zarrin, Mohammad, \emph{``Derived length and centralizers of groups"},  J. Algebra Appl. 14 (2015), no. 8, 1550133, 4 pp. 
\end{thebibliography}
\end{document}